\documentclass{amsart}

\usepackage{fancyhdr}
\usepackage{lastpage}
\usepackage{stmaryrd,yhmath}

\newcommand{\bdis}{\begin{displaymath}}
\newcommand{\edis}{\end{displaymath}}
\newcommand{\be}{\begin{equation}}
\newcommand{\ee}{\end{equation}}
\newcommand{\mbb}{\mathbb}
\newcommand{\mcal}{\mathcal}

\newcommand{\vp}{\varphi}

\newcommand{\vth}{\vartheta}

\newcommand{\zf}{\zeta\left(\frac{1}{2}+it\right)}

\newcommand{\okT}{\overset{k}{T}}

\theoremstyle{definition}

\theoremstyle{remark}
\newtheorem{remark}[]{Remark}

\newtheorem*{mydef1}{{\bf Theorem}}

\newtheorem*{mydef41}{{\bf Corollary 1}}

\newtheorem*{mydef42}{{\bf Corollary 2}}

\newtheorem*{mydef43}{{\bf Corollary 3}}

\newtheorem*{mydef44}{{\bf Corollary 4}}

\newtheorem*{mydef45}{{\bf Corollary 5}}

\newtheorem*{mydef6}{{\bf Example}}

\numberwithin{equation}{section}

\begin{document}

\title{Jacob's ladders and multiplicative algebra of reversely iterated integrals (energies) on the critical line}

\author{Jan Moser}

\address{Department of Mathematical Analysis and Numerical Mathematics, Comenius University, Mlynska Dolina M105, 842 48 Bratislava, SLOVAKIA}

\email{jan.mozer@fmph.uniba.sk}


\begin{abstract}
Certain completely logarithmic formula for a set of reversely iterated integrals (energies) is proved in this paper.
Namely, in this case we have that integral powers of $\ln T$ are contained on input as well as on output of
corresponding integrals (energies).
\end{abstract}
\maketitle

\section{Introduction}

Let us define the following matrix of reversely iterated segments (comp. \cite{3}, (2.3))
\be
||[\okT,\overset{q}{\wideparen{T+\ln^{-p}T}}]||_{p,q},\quad p,q=1,\dots,k,\ k\leq k_0,
\ee
where
\be
\begin{split}
 & \vp_1\{ [\overset{q}{T},\overset{q}{\wideparen{T+\ln^{-p}T}}]\}=
 [\overset{q-1}{T},\overset{q-1}{\wideparen{T+\ln^{-p}T}}], \\
 & \overset{0}{T}=T,\overset{0}{\wideparen{T+\ln^{-p}T}}=T+\ln^{-p}T,
\end{split}
\ee
where $\vp_1(t)$ is the Jacob's ladder (see (\cite{1}, \cite{2}) and $k_0\in\mbb{N}$ is and arbitrary fixed number.

Next we asign  to the $p$th row of the matrix (1.1) the following disconnected set
\be
\Delta_p(T,\ln^{-p}T)=\bigcup_{q=1}^k [\overset{q}{T},\overset{q}{\wideparen{T+\ln^{-p}T}}],\ p=1,\dots,k.
\ee
Properties of these sets are listed below (comp. \cite{3}, (2.5) -- (2.7)): since
\bdis
\ln^{-p}T=o\left(\frac{T}{\ln T}\right),\ p=1,\dots,k
\edis
then
\be
|[\overset{q}{T},\overset{q}{\wideparen{T+\ln^{-p}T}}]|=
\overset{q}{\wideparen{T+\ln^{-p}T}}-\overset{q}{T}=o\left(\frac{T}{\ln T}\right), \
q=1,\dots,k,
\ee
\be
|[\overset{q-1}{\wideparen{T+\ln^{-p}T}},\overset{q}{T}]|\sim (1-c)\pi(T); \pi(T)\sim\frac{T}{\ln T},\ q=2,\dots,k,
\ee
\be
\begin{split}
 & [\overset{1}{T},\overset{1}{\wideparen{T+\ln^{-p}T}}]\prec [\overset{2}{T},\overset{2}{\wideparen{T+\ln^{-p}T}}]
 \prec \dots \prec \\
 & \prec [\overset{k}{T},\overset{k}{\wideparen{T+\ln^{-p}T}}],\ p=1,\dots,k.
\end{split}
\ee

\begin{remark}
 Consequently, the asymptotic behavior of each of disconnected sets (1.3) is as follows: if $T\to\infty$ then the
 components recede unboundedly each from other and all together are receding to infinity. Hence, each of the sets
 (1.3) behaves like an one-dimensional Friedmann-Hubble expanding universe.
\end{remark}

Now we define the following correspondence
\bdis
[\overset{q}{T},\overset{q}{\wideparen{T+\ln^{-p}T}}]\longrightarrow
\int_{\overset{q}{T}}^{\overset{q}{\wideparen{T+\ln^{-p}T}}}
\prod_{r=0}^{q-1}\left|\zeta\left( \frac 12+i\vp_1^r(t)\right)\right|^2{\rm d}t
\edis
acting on the set of elements of the matrix (1.1). In this paper we obtain a canonical formula for these
integrals.

\section{Theorem on exclusivity of integer powers of $\ln T$}

\subsection{}

Let us remind that we have proved the following theorem (see \cite{3}, (2.1) -- (2.7)): for every
$L_2$-orthogonal system
\bdis
\{ f_n(t)\}_{n=1}^\infty,\ t\in [0,2l],\ l=o\left(\frac{T}{\ln T}\right),\ T\to\infty
\edis
there is a continuum set of $L_2$-orthogonal systems
\be
\begin{split}
 & \{ F_n(t;T,k,l)\}_{n=1}^\infty= \\
 & = \left\{ f_n(\vp_1^k(t)-T)\prod_{r=0}^{k-1}|\tilde{Z}[\vp_1^r(t)]|\right\},\
 t\in [\overset{k}{T},\overset{k}{\wideparen{T+2l}}],
\end{split}
\ee
where
\bdis
\begin{split}
 & \vp_1\{ [\overset{k}{T},\overset{k}{\wideparen{T+2l}}]\}=
 [\overset{k-1}{T},\overset{k-1}{\wideparen{T+2l}}],\ k=1,\dots,k_0, \\
 & [\overset{0}{T},\overset{0}{\wideparen{T+2l}}]=[T,T+2l],\ T\to\infty,
\end{split}
\edis
i. e. the following formula is valid
\be
\begin{split}
 & \int_{\overset{k}{T}}^{\overset{k}{\wideparen{T+2l}}}
 f_m(\vp_1^k(t)-T)f_n(\vp_1^k(t)-T)\prod_{r=0}^{k-1}
 \tilde{Z}^2[\vp_1^r(t)]{\rm d}t= \\
 & = \left\{ \begin{array}{rcl} 0 & , & m\not=n , \\ A_n & , & m=n, \end{array} \right. \\
 & A_n=\int_0^{2l} f_n^2(t){\rm d}t.
\end{split}
\ee
Of course, we have that
\be
\begin{split}
& \tilde{Z}^2(t)=\frac{{\rm d}\vp_1(t)}{{\rm d}t}=\frac{Z^2(t)}{2\Phi'_\vp[\vp(t)]}=\frac{\left|\zf\right|^2}{\omega(t)},\
 \vp_1(t)=\frac 12\vp(t), \\
& \omega(t)=\left\{ 1+\mcal{O}\left(\frac{\ln\ln t}{\ln t}\right)\right\}\ln t.
\end{split}
\ee
and
\bdis
\begin{split}
 & Z(t)=e^{i\vth(t)}\zf, \\
 & \vth(t)=-\frac t2\ln \pi +\text{Im}\ln\Gamma\left(\frac 14+i\frac t2\right).
\end{split}
\edis
Hence, for the classical Fourier's orthogonal system
\be
\left\{ 1,\cos\frac{\pi t}{l},\sin\frac{\pi t}{l},\dots,\cos\frac{n\pi t}{l},\sin\frac{n\pi t}{l},\dots\right\},\ t\in [0,2l],
\ee
for example, we have as the corresponding continuum set of orthogonal systems the following (see (2.1))
\be
\begin{split}
& \left\{ \prod_{r=0}^{k-1} \frac{\left|\zeta\left(\frac 12+i\vp_1^r(t)\right)\right|}{\sqrt{\omega[\vp_1^r(t)]}}\right.,\dots, \\
& \left. \left(\prod_{r=0}^{k-1} \frac{\left|\zeta\left(\frac 12+i\vp_1^r(t)\right)\right|}{\sqrt{\omega[\vp_1^r(t)]}}\right)
\cos\left(\frac{\pi}{l}n(\vp_1^k(t)-T)\right), \right. \\
& \left. \left(\prod_{r=0}^{k-1} \frac{\left|\zeta\left(\frac 12+i\vp_1^r(t)\right)\right|}{\sqrt{\omega[\vp_1^r(t)]}}\right)
\cos\left(\frac{\pi}{l}n(\vp_1^k(t)-T)\right),\dots \right\}, \\
& t\in [\overset{k}{T},\overset{k}{\wideparen{T+2l}}], \ k=1,\dots,k_0, \\
& T\in [T[\vp_1],+\infty).
\end{split}
\ee

\subsection{}

We have already noticed in our paper \cite{3} that the formula (2.2) can serve as a resource for new integral identities in the theory of the
Riemann zeta-function. We will see that this actually holds true. Either in the simplest case of the first function of the Fourier's system (2.4)
\bdis
f_1(t)=1
\edis
we have the following formula
\bdis
\int_{\overset{k}{T}}^{\overset{k}{\wideparen{T+2l}}}\prod_{r=0}^{k-1}\tilde{Z}^2[\vp_1^r(t)]{\rm d}t=2l,
\edis
i.e. (see (2.3),(2.5))
\be
\int_{\overset{k}{T}}^{\overset{k}{\wideparen{T+2l}}}
\prod_{r=0}^{k-1}\frac{\left|\zeta\left( \frac 12+i\vp_1^r(t)\right)\right|^2}{\omega[\vp_1^r(t)]}{\rm d}t=2l.
\ee
Next, let us remind that the disconnected set (see \cite{3}, (2.9))
\bdis
\Delta(T,k,l)=\bigcup_{r=0}^k [\overset{r}{T},\overset{r}{\wideparen{T+2l}}]
\edis
has the following properties (see \cite{3}, (2.7), (2.10), (4.3))
\be
\begin{split}
& [T,T+2l]\prec [\overset{1}{T},\overset{1}{\wideparen{T+2l}}]\prec \dots \prec
[\overset{k}{T},\overset{k}{\wideparen{T+2l}}]\prec \dots , \\
& \vp_1^r(t)\in [\overset{k-r}{T},\overset{k-r}{\wideparen{T+2l}}],\ r=0,1,\dots,k, \\
& \ln t\sim \ln T,\ \forall\ t\in (T,\overset{k}{\wideparen{T+2l}}),\ k=1,\dots,k_0.
\end{split}
\ee
Now we obtain from (2.6) by the mean-value theorem and (2.7) that
\be
\int_{\overset{k}{T}}^{\overset{k}{\wideparen{T+2l}}}\prod_{r=0}^{k-1}\left|\zeta\left( \frac 12+i\vp_1^r(t)\right)\right|^2
{\rm d}t\sim 2l\ln^kT,\ T\to\infty .
\ee
Consequently, we obtain from (2.8) in the case
\bdis
2l=\ln^{-p}T=o\left(\frac{T}{\ln T}\right),\ k=q
\edis
the following formula (comp. (1.7))
\begin{mydef1}
\be
\begin{split}
& \int_{\overset{q}{T}}^{\overset{q}{\wideparen{T+\ln^{-p}T}}}\prod_{r=0}^{q-1}\left|\zeta\left( \frac 12+i\vp_1^r(t)\right)\right|^2
{\rm d}t\sim \ln^{q-p}T, \\
& p,q=1,\dots,k,\ k\leq k_0,\ T\to\infty,
\end{split}
\ee
where
\be
\begin{split}
& \vp_1^r(t):\ \vp_1^0(t)=t,\ \vp_1^1(t)=\vp_1(t),\ \vp_1^2(t)=\vp_1(\vp_1(t)), \dots \\
& \vp_1^r(t)\in [\overset{q-r}{T},\overset{q-r}{\wideparen{T+\ln^{-p}T}}], \ r=0,1,\dots ,q
\end{split}
\ee
(see (2.7)) and $k_0\in\mbb{N}$ is an arbitrary and fixed number.
\end{mydef1}

\begin{remark}
The formula (2.9) is the first completely logarithmic formula in the theory of the Riemann zeta-function in the following sense
\bdis
[\overset{q}{T},\overset{q}{\wideparen{T+\ln^{-p}T}}]\overset{(2.9)}{\rightarrow} \{ 1+o(1)\}\ln^{q-p}T,\ T\to\infty.
\edis
Namely, the integer powers of $\ln T$ are contained on input as well as on the output of the integral (2.9). Of course, the formula
(2.9) is not accessible by the current methods in the theory of the Riemann zeta-function.
\end{remark}

\section{Interpretations of the iterated integrals under study}

\subsection{}

We define the following planar figures
\be
\begin{split}
& S^k_{p,q}(T)=\left\{ (t,y):\ t\in [\overset{q}{T},\overset{q}{\wideparen{T+\ln^{-p}T}}],  \right. \\
& \left. y\in \left[ 0,\prod_{r=0}^{q-1}\left|\zeta\left( \frac 12+i\vp_1^r(t)\right)\right|^2\right]\right\}, \\
& p,q=1,\dots,k,\ T\to\infty.
\end{split}
\ee
Then
\be
\int_{\overset{q}{T}}^{\overset{q}{\wideparen{T+\ln^{-p}T}}}\prod_{r=0}^{q-1}\left|\zeta\left( \frac 12+i\vp_1^r(t)\right)\right|^2
{\rm d}t=m\{ S^k_{p,q}(T)\}.
\ee

\begin{remark}
By (3.2) we have usual geometric interpretation of the integrals in the formula (3.2) as the measures of corresponding
planar figures (3.1), and, of course, (see (2.9), (3.2))
\be
m\{ S^k_{p,q}(T)\}\sim \ln^{q-p}T,\quad T\to\infty.
\ee
\end{remark}

\subsection{}

Let us consider an oscillating process (of any nature) described by the function
\be
f_{p,q}^k(t)=\left\{
\begin{array}{rcl} \prod_{r=0}^{q-1}\left|\zeta\left( \frac 12+i\vp_1^r(t)\right)\right| & , & t\in
[\overset{q}{T},\overset{q}{\wideparen{T+\ln^{-p}T}}] \\
0 & , & \mbox{otherwise} . \end{array}\right.
\ee
Then we have by means of the Plancherel's $L_2$-theory that
\be
\begin{split}
& \int_0^\infty \{ f_{p,q}^k(t)\}^2{\rm d}t=
\int_{\overset{q}{T}}^{\overset{q}{\wideparen{T+\ln^{-p}T}}}\prod_{r=0}^{q-1}\left|\zeta\left( \frac 12+i\vp_1^r(t)\right)\right|^2{\rm d}t=
\int_0^\infty \{ F_{p,q}^k(\omega)\}^2{\rm d}\omega,
\end{split}
\ee
where
\bdis
F_{p,q}^k(\omega)=\sqrt{\frac{2}{\pi}}\int_0^\infty f_{p,q}^k(t)\cos(\omega t){\rm d}t
\edis
is the Fourier's cosine transformation of $f_{p,q}^k(t)$ and
\bdis
\{ F_{p,q}^k(\omega)\}^2{\rm d}\omega
\edis
is the energy corresponding to the interval of frequencies
\bdis
[\omega,\omega+{\rm d}\omega).
\edis

\begin{remark}
Now we have by (3.5) the following energetic interpretation
\be
\int_{\overset{q}{T}}^{\overset{q}{\wideparen{T+\ln^{-p}T}}}\prod_{r=0}^{q-1}\left|\zeta\left( \frac 12+i\vp_1^r(t)\right)\right|^2{\rm d}t=
E^k_{p,q}(T) ,
\ee
where $E^k_{p,q}(T)$ is the total energy of the oscillating process (3.4) and (see (2.9), (3.6))
\be
E^k_{p,q}(T)\sim \ln^{q-p}T,\quad T\to\infty.
\ee
\end{remark}

\section{On generators of the main set of energies and constraints on behavior of $\zf$}

Since (see (2.9))
\bdis
\begin{split}
& \int_{\overset{q}{T}}^{\overset{q}{\wideparen{T+\ln^{-p}T}}}\prod_{r=0}^{q-1}\left|\zeta\left( \frac 12+i\vp_1^r(t)\right)\right|^2{\rm d}t
\sim \ln^{q-p}T, \\
& \int_{\overset{Q}{T}}^{\overset{Q}{\wideparen{T+\ln^{-P}T}}}\prod_{r=0}^{Q-1}\left|\zeta\left( \frac 12+i\vp_1^r(t)\right)\right|^2{\rm d}t\sim
\ln^{Q-P}T
\end{split}
\edis
we obtain  the following

\begin{mydef41}
For every fixed
\bdis
(P,Q):\ P\not=Q,\quad P,Q=1,\dots,k
\edis
we have that
\be
\begin{split}
& \int_{\overset{q}{T}}^{\overset{q}{\wideparen{T+\ln^{-p}T}}}\prod_{r=0}^{q-1}\left|\zeta\left( \frac 12+i\vp_1^r(t)\right)\right|^2{\rm d}t \sim \\
& \sim
\left\{
\int_{\overset{Q}{T}}^{\overset{Q}{\wideparen{T+\ln^{-P}T}}}\prod_{r=0}^{Q-1}\left|\zeta\left( \frac 12+i\vp_1^r(t)\right)\right|^2
{\rm d}t\right\}^{\frac{q-p}{Q-P}}, \\
& p,q=1,.\dots,k,\quad T\to\infty.
\end{split}
\ee
\end{mydef41}

\begin{remark}
We have by (4.1) that every fixed energy with $P\not=Q$ is the generator of all over main set of energies that correspond to
$p,q=1,\dots,k$.
\end{remark}

The energies in (4.1) correspond to the segments
\be
[\overset{q}{T},\overset{q}{\wideparen{T+\ln^{-p}T}}],
[\overset{Q}{T},\overset{Q}{\wideparen{T+\ln^{-P}T}}].
\ee

\begin{remark}
If $q\not=Q$ then we see that big distance (comp. (1.4), (1.5))
\be
> A\frac{T}{\ln T}\to\infty,\quad T\to\infty
\ee
separates the segments (4.2).
\end{remark}

\begin{remark}
Now, we will give an interpretation of the set of formulae (4.1) as a continuum set of constraints on behavior of the Riemann function
\bdis
\zf,\quad t\to\infty.
\edis
In this direction we have that by the constraints (4.1) is expressed a hight degree of inner binding of the set of values of the
function $\zf$, although at big distances (see (4.3)).
\end{remark}

Next, if we use the mean-value theorem in (4.1) then we obtain the following.

\begin{mydef42}
There are numbers
\bdis
\begin{split}
& d_1^r=\vp_1^r(c_1)\in (\overset{q-r}{T},\overset{q-r}{\wideparen{T+\ln^{-p}T}}),\ r=0,1,\dots,q-1, \\
& d_2^r=\vp_1^r(c_2)\in (\overset{Q-r}{T},\overset{Q-r}{\wideparen{T+\ln^{-p}T}}),\ r=0,1,\dots,Q-1,
\end{split}
\edis
(comp. (2.10)), of course,
\bdis
d_1^0=\vp_1^0(c_1)=c_1,\ d_2^0=\vp_1^0(c_2)=c_2,
\edis
such that
\be
\begin{split}
& |[\overset{q}{T},\overset{q}{\wideparen{T+\ln^{-p}T}}]|\prod_{r=0}^{q-1}\left|\zeta\left(\frac 12+id_1^r\right)\right|^2\sim \\
& \sim \left\{|[\overset{Q}{T},\overset{Q}{\wideparen{T+\ln^{-P}T}}]|\right\}^{\frac{q-p}{Q-P}}
\prod_{r=0}^{Q-1}\left|\zeta\left(\frac 12+id_2^r\right)\right|^{2\frac{q-p}{Q-P}}, \ T\to\infty.
\end{split}
\ee
\end{mydef42}

\begin{remark}
The formula (4.4) expresses one kind of continuum set of constraints that we have mentioned above. Namely, certain type of constraints corresponds
to each formula that follows from (4.1).
\end{remark}

\section{Law of multiplication of energies}

Next, we obtain from (4.1) the following

\begin{mydef43}
\be
\begin{split}
& \int_{\overset{q}{T}}^{\overset{q_1}{\wideparen{T+\ln^{-p_1}T}}}\prod_{r=0}^{q_1-1}\left|\zeta\left( \frac 12+i\vp_1^r(t)\right)\right|^2{\rm d}t
\times \\
& \times \int_{\overset{q_2}{T}}^{\overset{q_2}{\wideparen{T+\ln^{-p_2}T}}}\prod_{r=0}^{q_2-1}\left|\zeta\left( \frac 12+i\vp_1^r(t)\right)\right|^2{\rm d}t
\sim \\
& \sim
\left\{\int_{\overset{Q}{T}}^{\overset{P}{\wideparen{T+\ln^{-P}T}}}\prod_{r=0}^{Q-1}
\left|\zeta\left( \frac 12+i\vp_1^r(t)\right)\right|^2{\rm d}t\right\}^{\frac{q_1+q_2-(p_1+p_2)}{Q-P}},\ T\to\infty,
\end{split}
\ee
where
\be
\begin{split}
& -k+1\leq q_1+q_2-(p_1+p_2)\leq k-1, \\
& p_1,q_1,p_2,q_2=1,\dots,k.
\end{split}
\ee
\end{mydef43}

\begin{remark}
We see that the main set of energies is not closed with respect to multiplication of the kind (5.1). Namely, if
\bdis
(p_1,q_1)=(1,k),\ (p_2,q_2)=(k,1),
\edis
for example, then (see (5.2))
\bdis
q_1+q_2-(p_1+p_2)=2k-2>k-1,\ k\geq 2.
\edis
However, in this case of constraints it follows from (5.1), (comp. (4.4)) that we may use the condition
\bdis
-k_0+1\leq q_1+q_2-(p_1+p_2)\leq k_0-1
\edis
instead of (5.2).
\end{remark}

\section{Unit energies and non-local equivalences of them}

Next, we obtain from (2.9) or (4.1) at $p=q$ the following

\begin{mydef44}
\be
\int_{\overset{p}{T}}^{\overset{p}{\wideparen{T+\ln^{-p_1}T}}}\prod_{r=0}^{p-1}\left|\zeta\left( \frac 12+i\vp_1^r(t)\right)\right|^2{\rm d}t
\sim 1,\quad p=1,\dots,k,\ T\to\infty ,
\ee
i. e. these integrals play a role of the asymptotic  unit elements in the main set of energies.
\end{mydef44}

\begin{remark}
Since (see \cite{3}, (5.6))
\be
\overset{k}{T}-\overset{k-1}{T}\sim (1-c)\frac{T}{\ln T},\ T\to\infty, k=1,\dots,k_0 ,
\ee
then we obtain from (6.1) the following non-local equivalences of the unit energies
\be
\begin{split}
& \int_{\overset{1}{T}}^{\overset{1}{\wideparen{T+\ln^{-1}T}}}\left|\zeta\left( \frac 12+it\right)\right|^2{\rm d}t\sim
\int_{\overset{2}{T}}^{\overset{2}{\wideparen{T+\ln^{-2}T}}}\prod_{r=0}^{1}\left|\zeta\left( \frac 12+i\vp_1^r(t)\right)\right|^2{\rm d}t\sim \\
& \sim \dots \sim
\int_{\overset{k}{T}}^{\overset{k}{\wideparen{T+\ln^{-k}T}}}\prod_{r=0}^{k-1}\left|\zeta\left( \frac 12+i\vp_1^r(t)\right)\right|^2{\rm d}t
\end{split}
\ee
together with corresponding set of constraints (comp. (4.4)).
\end{remark}

\begin{remark}
We see that the initial unit energy transmission (see (6.3))
\be
\int_{\overset{1}{T}}^{\overset{1}{\wideparen{T+\ln^{-1}T}}}\left|\zeta\left( \frac 12+it\right)\right|^2{\rm d}t\sim 1.
\ee
corresponds to the translations on big distances (comp. (6.2)) 
\bdis
\overset{1}{T}\rightarrow\overset{2}{T}\rightarrow\dots\rightarrow\overset{k}{T}
\edis
Thus, the operation
\bdis
\int_{\overset{p}{T}}^{\overset{p}{\wideparen{T+\ln^{-p}T}}}\prod_{r=0}^{p-1}
\edis
asymptotically preserves the value of unit energy.
\end{remark}

\section{Inverse energies}

Next, we obtain from (5.1) in the case
\bdis
(p_1,q_1)=(p,q),\ (p_2,q_2)=(q,p)
\edis
the following

\begin{mydef45}
\be
\begin{split}
& \int_{\overset{q}{T}}^{\overset{q}{\wideparen{T+\ln^{-p}T}}}\prod_{r=0}^{q-1}\left|\zeta\left( \frac 12+i\vp_1^r(t)\right)\right|^2{\rm d}t\times \\
& \times
\int_{\overset{p}{T}}^{\overset{p}{\wideparen{T+\ln^{-q}T}}}\prod_{r=0}^{p-1}\left|\zeta\left( \frac 12+i\vp_1^r(t)\right)\right|^2{\rm d}t\sim 1,\quad
p,q=1,\dots,k,\ T\to\infty,
\end{split}
\ee
i.e. the inverse energy to
\bdis
\int_{\overset{q}{T}}^{\overset{q}{\wideparen{T+\ln^{-p}T}}}\prod_{r=0}^{q-1}\left|\zeta\left( \frac 12+i\vp_1^r(t)\right)\right|^2{\rm d}t
\edis
is the following one
\bdis
\int_{\overset{p}{T}}^{\overset{p}{\wideparen{T+\ln^{-q}T}}}\prod_{r=0}^{p-1}\left|\zeta\left( \frac 12+i\vp_1^r(t)\right)\right|^2{\rm d}t ,
\edis
and vice versa.
\end{mydef45}

\begin{mydef6}
\bdis
\begin{split}
& \left\{
\int_{\overset{257}{T}}^{\overset{257}{\wideparen{T+\ln^{-1}T}}}\prod_{r=0}^{256}\left|\zeta\left( \frac 12+i\vp_1^r(t)\right)\right|^2{\rm d}t
\right\}^{-1}\sim \\
& \sim
\int_{\overset{1}{T}}^{\overset{1}{\wideparen{T+\ln^{-257}T}}}\prod_{r=0}^{256}\left|\zeta\left( \frac 12+i\vp_1^r(t)\right)\right|^2{\rm d}t,\
T\to\infty.
\end{split}
\edis
\end{mydef6}

\begin{remark}
We see that each energy in (7.1) is asymptotically balanced by another one (a kind of balance on a lever). Moreover, let us remind that there is a
system of constraints corresponding to (7.1) (comp. (4.4)).
\end{remark}

\thanks{I would like to thank Michal Demetrian for his help with electronic version of this paper.}

\end{document}